\documentclass[12pt]{article}
\usepackage{stmaryrd}
\usepackage{mathrsfs}
\usepackage{amsmath}
\usepackage{amsmath,amsthm,amssymb}
\usepackage{amsfonts}
\usepackage{latexsym}
\date{}

\date{}

\begin{document}
\title{Ringel-Hall Algebras of Duplicated Tame Hereditary Algebras$^\star$}

\author{{Hongchang Dong,   Shunhua Zhang}\\
        {\small School of Mathematics, Shandong University,
        Jinan, 250100, P.R.China}}
\maketitle
\begin{center}
\begin{minipage}{120mm}

{\bf Abstract. } { Let $A$ be a tame hereditary algebra over a
finite field $k$ with $q$ elements, and ${\overline{A}}$ be the
duplicated algebra of $A$. In this paper, we investigate the
structure of Ringel-Hall algebra $\mathscr{H} (\overline{A})$ and of
the corresponding composition algebra $\mathscr{C} (\overline{A})$.
As an application, we prove the existence of Hall polynomials
$g_{XY}^M$ for any $\overline{A}$-modules $M, X$ and $Y$ with $X$
and $Y$ indecomposable if $A$ is a tame quiver $k$-algebra, then we
also obtain some Lie subalgebras induced by $\overline{A}$.}

\vskip 0.2in


\end{minipage}
\end{center}


\footnote{ $^\star$Supported by the NSF of China (Grant No.
10771112) and NSF of Shandong Province (Grant No. Y2008A05).}

\footnote {Email\ addresses: hcdong@mail.sdu.edu.cn(H.Dong),
shzhang@sdu.edu.cn(S.Zhang)}

\section {Introduction}

The duplicated algebras are interesting algebras that have been introduced recently
in the context of cluster categories. In particular, it is the interesting theory of these algebras
which is relevant for the connection to cluster categories, we refer to [ABST1] for details.

\vskip 0.2in

Ringel-Hall algebras of finitary rings were introduced by Ringel [R2, R3] in order to deal with
possible filtrations of modules with fixed factors. It turns out that Ringel-Hall algebra approach
provides a nice framework for the realization of quantized enveloping algebras and Kac-Moody
algebras, see [Lu, G, PX, R3-R5]. Later, some fundamental structures for Ringel-Hall algebras
of hereditary algebras were proved, see, e.g., [PZ1, PZ2, ZZ, SZ1, SZ2].

\vskip 0.2in

In this paper, we investigate the structure of Ringel-Hall algebras
of duplicated tame hereditary algebras in section 3, and in section 4 we prove the
existence of Hall polynomials $g_{XY}^M$ for any $\overline{A}$-modules $M, X$ and $Y$
with $X$ and $Y$ indecomposable (Theorem 4.8) when $A$ is a tame quiver algebra. As an application,
in section 5 we also obtain some Lie subalgebras induced by duplicated tame quiver
algebras.  Section 2 is devoted to some notations and definitions needed for our research.

\vskip 0.2in

\section {Preliminaries}

\vskip 0.2in

Let $A$ be a finite dimensional algebra over a field $k$. We
denote by $A$-mod the category of finitely generated left
$A$-modules, and $A$-ind a full subcategory of $A$-mod containing
exactly one representative of each isomorphism class of
indecomposable $A$-modules. Given a class $\cal C$
of $A$-modules, we denote by ${\rm add}\ \cal C$
the subcategory of $A$-mod whose objects are
the direct summands of finite direct sums of modules in $\cal C$.
We denote by $\Gamma_A$ the Auslander-Reiten quiver
of $A$ and by $\tau$ the Auslander-Reiten translation of $A$.
We refer to [ARS, DR, R1] for further notations and definitions in representation theory.

\vskip 0.2in

Let $M$ and $N$ be indecomposable $A$-modules. A path from $M$ to $N$ in $A$-ind
 is a sequence of non-zero morphisms
  $$M=M_0\stackrel{f_1} \longrightarrow M_1\stackrel{f_2} \longrightarrow\cdots
  \stackrel{f_t} \longrightarrow M_t =N$$
  with all $M_i$ in $A$-ind. Following [R1], we denote the
  existence of such a path by $M\leq N$. We say that $M$ is a
  predecessor of $N$ (or that $N$ is a successor of $M$).

\vskip 0.2in

More generally, if $S_1$ and $S_2$ are two sets of modules, we
 write $S_1\leq S_2$ if every module in $S_2$ has a predecessor in
 $S_1$, every module in $S_1$ has a successor in
 $S_2$, no  module in $S_2$ has a successor in
 $S_1$ and no  module in $S_1$ has a predecessor in
 $S_2$. The notation $S_1<S_2$ stands for $S_1\leq S_2$
 and $S_1\cap S_2 = \emptyset.$

\vskip 0.2in

Given a finite set $M$, we denote its cardinality by $|M|$. In the
sequel, we always assume that $k$ is a finite field with $q$
elements, that is $|k|=q$, and assume that $A$ is a
finite-dimensional tame hereditary algebra over $k$.

\vskip 0.2in

Let $M, N_1, N_2$ be finite dimensional $A$-modules. We denote by
$G_{N_1N_2}^{M}$ the number of submodules $L$ of $M$ with the
property that $L\simeq N_2$ and $M/L \simeq N_1$. By [R2] the
Ringel-Hall algebra ${\mathscr H}(A)$ is a free abelian group with a
basis ${\{u_{[M]}}\}_{[M]}$ indexed by the isomorphism classes of
finite (left) $A$-modules with the multiplication defined by
$$
u_{[N_1]}\cdot u_{[N_2]} = \sum\limits_{[M]} G_{N_1N_2}^{M} u_{[M]}
$$
Note that we only deal with finite sum since $A$ is a finite ring.
We denote by ${\mathscr C}(A)$ the subalgebra of ${\mathscr H}(A)$
generated by simple $A$-modules which is called composition algebra.

\vskip 0.2in

From now on, we always assume that $A$ is a tame hereditary $k$-algebra
and that $\overline{A}$ is the duplicated algebra of $A$, see [ABST1].
Then $\overline{A}= ( \begin{array}{lr} A & 0\\
DA & A \end{array} )$ is the matrix algebra, we see that
$\overline{A}$ contains two copies of $A$ given respectively by
$e\overline{A}e$  and by $e'\overline{A}e'$,  where $e=( \begin{array}{lr} 1 & 0\\
0 & 0 \end{array} )$,  and $e'=( \begin{array}{lr} 0 & 0\\
0 & 1 \end{array} )$.  We denote the first one by $A=A_0$ and the second
one by $A'=A_1$.  Accordingly,  $Q'_A$ denotes the quiver of $A'$, $x'$
the vertex of $Q'_A$ cerresponding to $x\in (Q_A)_0$,  and $e'_x$
the corresponding idempotent.  Let $S_x, P_x, I_x$ denote respectively the corresponding simple,
 indecomposable projective and indecomposable injective module in
 $\overline{A}$-mod  corresponding to $x\in(Q_A\bigcup Q'_A)_0$. For simpleness, we write
 $Q_A=\{ 1,\cdots, n \}$ and $Q'_A=\{ 1',\cdots, n' \}$.

\vskip 0.2in

Recall from [ABST2] that the Auslander-Reiten quiver of
$\Gamma_{\overline{A}}$ can be described as follows. It starts with
the Auslander-Reiten quiver of $A_0=A$. Then projective-injective
modules start to appear, such projective-injective module has its
socle corresponding to a simple $A_0$-module, and its top
corresponding to a simple $A_1$-module. Next occurs a part denoted
by $A_{01}$-ind where indecomposables contain at same time simple
composition factors from simple $A_0$-modules, and simple
composition factors from simple $A_1$-modules. When all
projective-injective modules whose socle corresponding to simple
$A_0$-modules have appeared, we reach the projective $A_1$-modules
and thus the Auslander-Reiten quiver of $A_1$.

\vskip0.2in

From the description above, we can divide the
Auslander-Reiten-quiver $\Gamma_{\overline{A}}$ into 7 parts,
denoted by $\mathscr{P}_0$, $\mathscr{R}_0$, $\mathscr{X}_0$,
$\mathscr{R}_{01}$, $\mathscr{X}_1$, $\mathscr{R}_1$,
$\mathscr{I}_1$ respectively, where $\mathscr{P}_0$ (resp.
$\mathscr{I}_1$)is the preprojective (resp. preinjective) component
of $\Gamma_{A_0}$ (resp. $\Gamma_{A_1}$); $\mathscr{X}_0$ and
$\mathscr{X}_1$ are forms of translation quiver of $\mathbb{Z}
Q_{A}$; $\mathscr{R}_0$, $\mathscr{R}_{01}$ and $\mathscr{R}_1$ are
the same types of tubes since $A$ is tame type.

\vskip 0.2in

Let $M$ be an $\overline{A}$-module. We denote by $\Omega^{-i}(M)$
the $i^{\rm th}$ cosyzygy of $M$ and by $\Omega^{i}(M)$
the $i^{\rm th}$ syzygy of $M$ respectively.
Let ${\cal L}_{\overline{A}}$ be the left part of $\overline{A}$-mod.
By definitions in [HRS], ${\cal L}_{\overline{A}}$ is the full subcategory
of ${\overline{A}}$-mod consisting of all indecomposable ${\overline{A}}$-modules
such that if $L$ is a predecessor of $M$, then the projective dimension ${\rm pd}\ L$
of $L$ is at most one.

\vskip 0.2in

It is well known that ${\rm gl.dim} \overline{A}$, the global dimension of $\overline{A}$,
is 3.  We denote by $\Sigma_0$ the set of all non-isomorphic indecomposable
projective $A_0$-modules,  and write
 $\Sigma_k=\Omega^{-k}\,\Sigma_0=\{ \Omega^{-k} X \  | \ X\in \Sigma_0 \}$ for $1\leq k\leq 2$.

\vskip 0.2in

\section {Ringel-Hall algebras of duplicated tame hereditary algebras}

In this section, we mainly investigate the structure of Ringel-Hall
algebra of duplicated tame hereditary algebras, the decomposition of
the composition algebra, indecomposable elements in the composition
algebra, and prove that the exceptional elements can be written as
skew communicators.

\vskip0.2in

Let $A_{01}$-ind be the indecomposable $\overline{A}$-modules with
the composition factors having both simple $A_0$-modules and
$A_1$-modules. That is $\overline{A}$-ind = $A_0$-ind $\bigcup
A_{01}$-ind $\bigcup A_1$-ind. Note that ${\rm add} (A_{01}$-ind) is
an exact category which is closed under extensions, we can define
the corresponding Ringel-Hall algebra which is denoted by
$\mathscr{H}(A_{01})$.

\vskip0.2in

{\bf Theorem 3.1.} \ {\it
 $\mathscr{H}(\overline{A})= \mathscr{H}(A_0) \mathscr{H}(A_{01}) \mathscr{H}(A_1)$}.

\vskip0.1in

 {\bf Proof.}  Let $M$ be an $\overline{A}$-module. Suppose
that $[M]= [M_0]\oplus [M_{01}]\oplus [M_1]$, where $M_{01}\in {\rm
add}\ (A_{01}-{\rm ind})$, and $M_i\in A_{i}-{\rm mod}$ with $i= 0, 1$.

It is easy to see that $g^{ M_{01}\oplus M_1}_{ M_{01}\  M_1}= 1$,
$g^ M_{ M_0\ M_{01}\oplus M_1}= 1$, and
$$
u_{[M_0]} u_{[M_{01}]}u_{[M_1]} = u_{[M_0]} u_{[M_{01}]\oplus
[M_1]}=u_{[M_0]\oplus [M_{01}]\oplus [M_1]}=u_{[M]}.
$$
                      $\hfill\Box$

\vskip0.2in

Let $\mathscr{C}(A_{01})= \mathscr{C}(\overline{A})\bigcap \mathscr{H}(A_{01})$
and $\mathscr{C}(A_{0})$ (resp. $\mathscr{C}(A_{1})$) be the subalgebra of
$\mathscr{C}(\overline{A})$ generated by the simple modules $S_1, \cdots, S_n$
(resp. $S_{1'}, \cdots, S_{n'}$). Note that
$\mathscr{C}(A_{0})= \mathscr{H}(A_0)\bigcap \mathscr{C}(\overline{A})$ and
that $\mathscr{C}(A_{1})= \mathscr{H}(A_1)\bigcap \mathscr{C}(\overline{A})$,
by using Theorem 3.1, we have the following.

\vskip0.2in

{\bf Corollary 3.2.} \ {\it $\mathscr{C}(\overline{A})=
\mathscr{C}(A_{0})\mathscr{C}(A_{01})\mathscr{C}(A_{1})$}.

\vskip0.2in

Let $M$ be an indecomposable $\overline{A}$-module.  $M$ is said to
be exceptional if ${\rm Ext}^i_{\overline{A}}(M,M)= 0$ for all
$i>0$.

\vskip0.2in

{\bf Theorem 3.3.} \ {\it Let  $\overline{A}$ be the duplicated tame
hereditary algebra $A$ and  $M$ be an indecomposable
$\overline{A}$-module. Then $u_{[M]}\in \mathscr{C}(\overline{A})$ if and
only if $M$ is an exceptional $\overline{A}$-module}.

\vskip0.1in

 {\bf Proof.}\ Assume that $M$ is an exceptional
$\overline{A}$-module.

\vskip0.1in

{\bf Case I.} If $M\in A_0-{\rm ind}$ or $M\in A_1-{\rm ind}$, then
$u_{[M]}\in \mathscr{C}(\overline{A})$ by [PZ1, SZ2].

\vskip0.1in

{\bf Case II.} Assume that $M\in A_{01}-{\rm ind}$. First of all, we
suppose that $M$ is a projective-injective $\overline{A}$-module.

If $M\in {\cal L}_{\overline{A}}$, then ${\rm top} M = S_{i'}$ and ${\rm rad}
M\in A_0-{\rm mod}$. Note that ${\rm rad} M$ is a preinjectve
$A_0$-module and $u_{[{\rm rad} M]}\in \mathscr{C}(\overline{A})$, we
have the following:
$$
u_{[M]} = u_{[S_{i'}]}u_{[{\rm rad} M]}-u_{[{\rm rad} M]}u_{[S_{i'}]}\in
\mathscr{C}(\overline{A}).
$$

If $M\not\in {\cal L}_A$, then ${\rm Soc} M = S_{i}$ and $M/{\rm
Soc} M\in A_1-{\rm mod}$. In this case, one can easy to see that
$M/{\rm Soc} M$ is a preprojectve $A_1$-module and $u_{[M/{\rm Soc}
M]}\in \mathscr{C}(\overline{A})$, we have the following:
$$
u_{[M]} = u_{[M/{\rm Soc} M]}u_{[S_{i}]}-u_{]S_{i}]}u_{]M/{\rm Soc} M]}\in
\mathscr{C}(\overline{A}).
$$

Finally, we can assume that $M\in A_{01}-{\rm ind}$ and $M$ is not a
projective-injective $\overline{A}$-module. Read from the
Auslander-Reiten quiver of $\overline{A}$ and by using Theorem 9.1
in [PZ1] and Theorem 1 in [SZ2], we know that $u_{[M]}\in
\mathscr{C}(\overline{A})$.

 \vskip0.1in

Conversely, let $M$ be an indecomposable $\overline{A}$-module and
$u_{[M]}\in \mathscr{C}(\overline{A})$. We want to prove that $M$ is
an exceptional $\overline{A}$-module. If $M\in A_0-{\rm ind}$ or
$M\in A_1-{\rm ind}$, then $M$ is an exceptional
$\overline{A}$-module follows from [ZZ] if $A$ is a tame quiver
algebra and follows from [SZ2] when $A$ is a non-simply-laced tame
hereditary algebra.

Now assume that $M\in A_{01}-{\rm ind}$ and we may assume that $M$
is not a projective-injective $\overline{A}$-module. It is easy to
read from the Auslander-Reiten quiver of $\overline{A}$ that $M$ is
in a full subquiver of $\Gamma_{\overline{A}}$ which is isomorphic
to $\Gamma_{D(A_0)}$. Then $M$ is an exceptional
$\overline{A}$-module follows from [ZZ, SZ2] again. This completes
the proof.  $\hfill\Box$

\vskip 0.2in

{\bf Example 3.4.}\  Let $\overline{A}$ be the duplicated tame
quiver algebra of type $\widetilde{D}_4$. That is,
$\overline{A}=\overline{k \widetilde{D_4}}/I$, and
$\overline{\widetilde{D_4}}$ is the following quiver,

$$\begin{array}{ccccccc}

       && 2 & & &&2'\\[-2ex]

       &\swarrow & &\nwarrow & &\swarrow& \\[-2ex]

  \overline{\widetilde{D_4}}:    1 &\leftleftarrows & \begin{array}{c}3\\[-2ex]4\end{array}
    &\leftleftarrows & 1' &\leftleftarrows &  \begin{array}{c}3'\\[-2ex]4'\end{array}\\[-2ex]

      &  \nwarrow & & \swarrow&&\nwarrow &\\[-2ex]

    && 5 &&&&5'\\[-2ex]

    \end{array}.
$$

\vskip 0.2in

Then the indecomposable projective-injective $\overline{A}$-modules
are represented by their Loewy series as the following,
$$
P_1'=\begin{array}{c}1'\\[-2ex]2345\\[-2ex]1\end{array},\
P_2'=\begin{array}{c}2'\\[-2ex]1'\\[-2ex]2\end{array},\
P_3'=\begin{array}{c}3'\\[-2ex]1'\\[-2ex]3\end{array}, \
P_4'=\begin{array}{c}4'\\[-2ex]1'\\[-2ex]4\end{array}, \
P_5'=\begin{array}{c}5'\\[-2ex]1'\\[-2ex]5\end{array}.
$$

We should mention that the minimal positive imaginary root of
$k\widetilde{D}_4$ is $\delta=(2,1,1,1,1)$ and every indecomposable
$\overline{A}$-module $M$ which belongs to $\mathscr{R}_0$,
$\mathscr{R}_{01}$ or $\mathscr{R}_1$ with $l(M)\geq 6$ is not
exceptional, where $l(M)$ is the length of $M$.

According to Theorem 3.3, for any indecomposable
$\overline{A}$-module $M$,  $u_{[M]}$ belongs to
$\mathscr{C}(\overline{A})$ if and only if $M$ belongs to
$\mathscr{P}_0$, $\mathscr{X}_0$, $\mathscr{X}_1$, $\mathscr{I}_1$
or to $\mathscr{R}_0$, $\mathscr{R}_{01}$, $\mathscr{R}_1$ with
$l(M)< 6$.

\vskip0.2in

The following concept is defined in [PZ2]. Let $B$ be a $k$-algebra,
and $x,y\in B$, and $c,d\in k^*= k\backslash \{ 0 \}$. The element
$cxy-dyx$ is called a {\it skew commutator} of $x$ and $y$. Let
$X=\{ x_1, \cdots, x_n\}$ be a set of $B$. Define the sets $X_i$
inductively: Let $X_0= X$. Let $X_i$ be the set of all skew
commutators of arbitrary two different elements in
$\bigcup\limits_{j<i} X_j$. An element $x\in B$ is called an
iterated skew commutator of $x_1, \cdots, x_n$, provided that there
exists a positive integer $m$ such that $x\in X_m$.

\vskip0.2in

The following theorem indicates that an indecomposable non-simple
$\overline{A}$-module which belongs to $ \mathscr{C}(\overline{A})$ can
be written as an iterated skew commutator of simple $\overline{A}$-modules.

\vskip0.2in

{\bf Theorem 3.5.} \ {\it Let $A$ be a tame hereditary algebra over
$k$ and $\overline{A}$ be the duplicated algebra of $A$. Let $M$ be
a non-simple indecomposable $\overline{A}$-module. Then the element
$u_{[M]}\in \mathscr{C}(\overline{A})$ can be written as an iterated
skew commutator of the isoclasses of simple $\overline{A}$-modules}.

\vskip0.1in

 {\bf Proof.}\  According to Theorem 3.3, we know that $M$ is an
 exceptional $\overline{A}$-module.

 If $M\in A_0-{\rm ind}$ or $M\in A_1-{\rm ind}$, then
$M$ is an iterated skew commutator of the isoclasses of simple
$\overline{A}$-modules by Theorem 2.1 in [PZ2].

Now, Let $M\in A_{01}-{\rm ind}$ and $M$ be a
projective-injective $\overline{A}$-module.

If $M\in {\cal L}_{\overline{A}}$, according to the proof
of Theorem 3.3, we can write $u_{[M]}$ as following:
$$
u_{[M]} = u_{[S_{i'}]}u_{[{\rm rad} M]}-u_{[{\rm rad} M]}u_{[S_{i'}]}\in
\mathscr{C}(\overline{A}),
$$
${\rm rad} M\in A_0-{\rm mod}$ is a preinjectve $A_0$-module and
${\rm top} M = S_{i'}$. By using the Theorem 2.1 in [PZ2], $u_{[{\rm
rad} M]}$ is an iterated skew commutator of the isoclasses of simple
$A_0$-modules, hence $u_{[M]}$ can be written as an iterated skew
commutator of the isoclasses of simple $\overline{A}$-modules.

If $M\not\in {\cal L}_{\overline{A}}$, then ${\rm Soc} M = S_{i}$ is
a simple  $A_0$-module and $M/{\rm Soc} M$  is a preprojectve
$A_1$-module.  By using Theorem 2.1 in [PZ2] again, $u_{[M/{\rm Soc}
M]}$ can be written as an iterated skew commutator of the $S_{1'},
\cdots, S_{n'}$.

Note that $u_{[M]} = u_{[M/{\rm Soc} M]}u_{[S_{i}]}-u_{[S_{i}]}u_{[M/{\rm Soc} M]}\in
\mathscr{C}(\overline{A})$, thus $u_{[M]}$ is an iterated skew commutator of the isoclasses of
simple $\overline{A}$-modules.

 Finally, we can assume that $M\in
A_{01}-{\rm ind}$ and $M$ is not a projective-injective
$\overline{A}$-module. It follows from the Auslander-Reiten quiver
of $\overline{A}$ that $M$ is in a full subquiver of
$\Gamma_{\overline{A}}$ which is isomorphic to $\Gamma_{D(A_0)}$.
Since $M$ is an exceptional $\overline{A}$-module,  we know that
$u_{[M]}$ is an iterated skew commutator of the isoclasses of simple
$\overline{A}$-modules by using Theorem 2.1 in [PZ2].  The proof is
completed.           $\hfill\Box$

\vskip0.2in

\section {Some Hall polynomials for duplicated tame hereditary algebras}

In this section, we always assume that $A$ is a tame quiver algebra
over $k$ and $\overline{A}$ be the duplicated algebra of $A$, and we
will prove that some Hall polynomials for duplicated tame hereditary
algebras exist. Note that we can, in this case, divide the
Auslander-Reiten quiver $\Gamma_{\overline{A}}$ into 7 parts,
denoted by $\mathscr{P}_0$, $\mathscr{R}_0$, $\mathscr{X}_0$,
$\mathscr{R}_{01}$, $\mathscr{X}_1$, $\mathscr{R}_1$,
$\mathscr{I}_1$ respectively.

\vskip 0.2in

Let $E$ be a field extension of $k$. For any
$k$-space $V$, we denote by $V^E$ the $E$-space $V\otimes_kE$; then,
of course, $A^E$ naturally becomes an $E$-algebra. If $S$ is a
simple $A$-module, according to Theorem 7.5 in [La], we know that
$S^E$ is the simple $A^E$-module. For any $M\in A$-mod, $E$ is
called $M$-conservative for $A$ if for any indecomposable summand
$N$ of $M$, $({\rm End}N/{\rm rad\ End}N)^E$ is a field. Under field
isomorphism, we put
$$\Omega_M=\{E|E \ {\rm is\ a\ finite\ field\
extension\ of}\ k\ {\rm and}\ E\ {\rm is}\ M-{\rm conservative\
for}\ A\}.$$

\vskip 0.2in

Note that $\Omega_M$ is an infinite set, since $M$ has only finitely
indecomposable summands. By [SZ1], we say that Hall polynomials exist
for $A$, if for any $M$, $N_1$, $N_2\in A$-mod, there exists a
polynomial $g_{N_1N_2}^M\in\mathbf{Z}[x]$ and an infinite subset
$\Omega_{N_1N_2}^M$ of $\Omega_{M\oplus N_1\oplus N_2}$, such that
for any $E\in\Omega_{N_1N_2}^M$,
$$g_{N_1N_2}^M(|E|)=G_{N_1^EN_2^E}^{M^E}.$$
Such a polynomial $g_{N_1N_2}^M$ is called a Hall polynomial of $A$.

\vskip 0.2in

{\bf Remark.}\  When $A$ is a representation-finite algebra, the
above definition is the same as in [R5].

\vskip0.2in

The following results were proved in [SZ1] for tame quiver algebras,
and we observe that they are also true for duplicated tame hereditary
algebras,  we refer to [SZ1] for details.

\vskip0.2in

{\bf Lemma 4.1.} \ {\it Assume that $u_{[M]}=u_{[M_1]}+u_{[M_2]}$ in
$\mathscr{H}(\overline{A})$. Then $G_{MN}^L=G_{M_1N}^L+G_{M_2N}^L$
for any $\overline{A}$-modules $L$ and $N$.}

\vskip0.2in

{\bf Lemma 4.2.} \ {\it
Given $M$, $N\in\overline{A}$-${\rm mod}$, then there exists a
nonnegative integer $h(M,N)$ such that $|{\rm
Hom}_{\overline{A}^E}(M^E,N^E)|=|E|^{h(M,N)}$ for any $E\in
\Omega_{M\oplus N}$.}

\vskip0.2in

{\bf Lemma 4.3.} \ {\it
Let $N_1,N_2,\cdots,N_t$ be simple $\overline{A}$-modules except at
most only one. Then there exists the Hall polynomial
$g_{N_1N_2\cdots N_t}^M$ for all $M\in\overline{A}$-${\rm mod}$.}

\vskip0.2in

{\bf Lemma 4.4.} \ {\it Let $M$, $N$, $L$ be $\overline{A}$-modules
with $N\in \mathscr{P}_0$, $\mathscr{X}_0$, $\mathscr{X}_1$ or $\mathscr{I}_1$.
Then the Hall polynomials $g_{NL}^M$ and  $g_{LN}^M$ exist.}

\vskip0.1in

{\bf Proof.} \  By duality, we only need to prove the existence of
the Hall polynomial $g_{NL}^M$. Note that for any indecomposable
$\overline{A}$-module $X\in \mathscr{P}_0$, $\mathscr{X}_0$,
$\mathscr{X}_1$ or $\mathscr{I}_1$ which is exceptional, according
to Theorem 3.3, we know that $u_{[X]}\in \mathscr{C}(\overline{A})$.
Therefore we can assume that
$$
u_{[N]}=\sum_{i_1,\cdots,i_t\in\{1,\cdots,n,1',\cdots,n'\}}a_{i_1\cdots
i_t}u_{[s_{i_1}]}\cdots u_{[s_{i_t}]},
$$
where $a_{i_1\cdots i_t}\in \mathbf{Z}$. By using Lemma 4.1, we have
$$
G_{NL}^{M}=\sum_{i_1,\cdots,i_t\in\{1,\cdots,n,1',\cdots,n'\}}a_{i_1\cdots
i_t}G_{S_{i_1}\cdots S_{i_t}L}^{M}.
$$
According to Lemma 4.3, we have Hall polynomials $g_{S_{i_1}\cdots S_{i_t}L}^M\in
\mathbf{Z}[x]$ such that there exists an infinite subset $\Omega_{NL}^M$ of $\Omega_{M\oplus
N\oplus L}$ and for any $E\in\Omega_{NL}^M$
$$
G_{S_{i_1}^E\cdots S_{i_t}^EL^E}^{M^E}=g_{S_{i_1}\cdots S_{i_t}L}^M(|E|).
$$
Let $g_{NL}^M=\sum_{i_1,\cdots,i_t\in\{1,\cdots,n,1',\cdots,n'\}}a_{i_1\cdots
i_t}g_{S_{i_1}\cdots S_{i_t}L}^{M}\in\mathbf{Z}[x]$.
For any $E\in\Omega_{NL}^M$, we have that
$G_{N^EL^E}^{M^E}=g_{NL}^M(|E|)$, that is, $g_{NL}^M$ is the Hall
polynomial of $\overline{A}$. This completes the proof.    $\hfill\Box$

\vskip0.2in

{\bf Lemma 4.5.} \ {\it Let $M$, $N$ and $L$ be $\overline{A}$-modules with $N$ and
$L$ indecomposable.  If  $N, L\in \mathscr{R}_{0}$, $N, L\in \mathscr{R}_{01}$ or
$N, L\in \mathscr{R}_{1}$, then the Hall polynomial $g_{NL}^M$ exists.}

\vskip0.1in

{\bf Proof.}\  We may assume that $M$ is an extension of $L$ by $N$,
since otherwise we may take $g_{NL}^M=0$. Therefore we have a short
exact sequence $0\rightarrow L\rightarrow M\rightarrow N\rightarrow
0$, it follows that $M, N, L$ belong to the same part of
$\Gamma_{\overline{A}}$, that is, $M, N, L\in \mathscr{R}_{0}$, $M,
N, L\in \mathscr{R}_{01}$ or  $M, N, L\in \mathscr{R}_{1}$. By using
the same method as Lemma 2.9 in [SZ1], we know that $g_{NL}^M$
exists.   This completes the proof.    $\hfill\Box$

\vskip0.2in
 For any $M, N, L \in \overline{A}-{\rm mod}$, we denote by
 ${\rm Ext}^1_{\overline{A}}(N,L)_M$ the set of all
exact sequences in ${\rm Ext}^1_{\overline{A}}(N,L)$ with middle term $M$.
The following lemma was proved in [P, Rie].

\vskip0.2in

{\bf Lemma 4.6.} \ {\it For any $\overline{A}$-modules $M$, $N$, $L$,
$$
G_{NL}^M=\frac{\mid{\rm
Ext}_{\overline{A}}^1(N,L)_M\mid\cdot\mid{\rm
Aut}_{\overline{A}}M\mid}{\mid{\rm Aut}_{\overline{A}}N\mid\cdot\mid
{\rm Aut}_{\overline{A}}L\mid\cdot\mid{\rm
Hom}_{\overline{A}}(N,L)\mid}.
$$}

\vskip0.2in

{\bf Lemma 4.7.} \ {\it Let $M, N$  and $L$ be
$\overline{A}$-modules with $N$ and $L$ indecomposable. Assume that
$N\in \mathscr{R}_{i}$, $L\in \mathscr{R}_{j}$ with  $i\neq
j\in\{0,01,1\}$. Then the Hall polynomial $g_{NL}^M$ exists.}

\vskip0.1in

{\bf Proof.} \  We only need to consider the cases $L\in \mathscr{R}_{0}$ with
 $N\in \mathscr{R}_{01}$, and $L\in \mathscr{R}_{01}$ with $N\in \mathscr{R}_{1}$,
 since in other cases we have that ${\rm Ext}_{\overline{A}}^1(N, L)=0$, and by using Lemma 4.6,
 the existence of the Hall polynomial $g_{NL}^M$ follows.

 \vskip0.1in

 {\bf Case I.}\ Let $L\in \mathscr{R}_{0}$ with $N\in \mathscr{R}_{01}$.
 Assume that $E(L)$ is the injective envelope of $L$, then we have a
 short exact sequence
$$
(*)\ \ \ \ \ \ \ \ \ \ \   0\rightarrow L\rightarrow E(L)\rightarrow \Omega^{-1}L\rightarrow 0,
$$
where $E(L)$ is projective-injective $\overline{A}$-module since $L\in \mathscr{R}_{0}$, and
$\Omega^{-1}L$ is an indecomposable $\overline{A}$-module which belongs to $\mathscr{R}_{01}$.
Note that $E(L)$ is a predecessor of $N$, by applying ${\rm Hom}_{\overline{A}}(N,-)$ to $(*)$,
we obtain that ${\rm Hom}_{\overline{A}}(N,\Omega^{-1}L)\simeq {\rm Ext}_{\overline{A}}^{1}(N,L)$.
Hence ${\rm dim}_k{\rm Ext}_{\overline{A}}^{1}(N,L)={\rm
dim}_k{\rm Hom}_{\overline{A}}(N,\Omega^{-1}L)\leq 1$ since $N$ and $\Omega^{-1}L$ are indecomposable
$\overline{A}$-modules belonging to $\mathscr{R}_{01}$.
For any $\overline{A}$-module $M$, according to Lemma 4.2 and Lemma 4.6 we know that the Hall
polynomial $g_{NL}^M$ exists.

 \vskip0.1in

 {\bf Case II.}\  Let $L\in \mathscr{R}_{01}$ with $N\in \mathscr{R}_{1}$.
 By using the same method as in Case I,  we can prove that the Hall
polynomial $g_{NL}^M$ exists.
This completes the proof.       $\hfill\Box$

\vskip0.2in

{\bf Theorem 4.8.} \ {\it Let $X$ and $Y$ be indecomposable
$\overline{A}$-modules. Then for any $\overline{A}$-module $M$,
there exists the Hall polynomial $g_{XY}^M$.}

\vskip0.1in

{\bf Proof.}\  If one of the indecomposable $\overline{A}$-modules $X$ and
$Y$ belongs to $\mathscr{P}_0$, $\mathscr{X}_0$, $\mathscr{X}_1$ or $\mathscr{I}_1$,
by Lemma 4.4, we know that the Hall polynomial $g_{NL}^M$ exists.

If none of $X$ and $Y$ belongs to $\mathscr{P}_0$, $\mathscr{X}_0$, $\mathscr{X}_1$ or $\mathscr{I}_1$,
then $X$ and $Y$ must belong to $\mathscr{R}_i$ or $\mathscr{R}_j$, where $(i,j\in\{0, 01,1\})$.
In case $i=j$, then the existence of the Hall polynomial $g_{XY}^M$ follows from Lemma 4.5.
If $i\neq j$, according to Lemma 4.7, we have the Hall polynomial $g_{XY}^M$ exists.
The proof is completed.                $\hfill\Box$

\vskip0.2in

{\bf Remark.}\ If $A$ is a representation-finite hereditary $k$-algebra, then the
duplicated algebra $\overline{A}$ is represented-direct, thus according to [R5],
we know that the Hall polynomial $g_{NL}^M$ exists for any $\overline{A}$-modules $M, N, L$.

\vskip0.2in

\section {The Lie subalgebras induced by duplicated tame hereditary algebras}

In this section, we also assume that $A$ is a tame quiver algebra over $k$
and $\overline{A}$ is the duplicated algebra of $A$, and we will investigate
some Lie subalgebras induced by $\overline{A}$ which seem to have an independent interest.

\vskip0.2in

Let $\Omega$ be an infinite set of finite field extension of $k$ up
to isomorphism. Since $A$ is a tame quiver algebra, according to
[CD] and  Theorem 7.5 in [La], we know that $E$ is $S$-conservative
for any simple $\overline{A}$-module $S$.

\vskip0.2in

Denote by ${\mathscr H}(\overline{A},\Omega)$ the
subring of $\prod\limits_{E\in\Omega}{\mathscr H}(\overline{A^E})$
generated by $\{([M^E])_{E\in\Omega}|M\in \overline{A}-{\rm mod}\}$
and $q_\Omega=(|E|u_{[0]})_{E\in\Omega}$. Denote by ${\mathscr
H}(\overline{A})_1$ the quotient ring ${\mathscr
H}(\overline{A},\Omega)/(q_\Omega-1){\mathscr
H}(\overline{A},\Omega)$, called the degenerate Ringel-Hall algebra
of $\overline{A}$. The subalgebra of ${\mathscr
H}(\overline{A})_1$ generated by the simple $\overline{A}$-modules,
denoted by $\mathscr{C}(\overline{A})_1$, is called the degenerate
composition algebra of $\overline{A}$.

\vskip0.2in

The following Lemma was proved in [R5].

\vskip0.2in

{\bf Lemma 5.1.}\ {\it Let $M, X, Y\in\overline{A}$-${\rm mod}$ with $X$ and
$Y$ indecomposable. For any $E\in\Omega_{M\oplus X\oplus Y}$, then

\vskip0.1in

${\rm(1)}$ If $M\not\simeq X\oplus Y$, then $|E|-1$ divides $G_{X^EY^E}^{M^E}$ $\rm;$

\vskip0.1in

${\rm(2)}$ If $M\simeq X\oplus Y$. If $X\simeq Y$, then $|E|-1$
divides $G_{X^EY^E}^{M^E}-2$;

\ \ \ \ If $X\not\simeq Y$, then $|E|-1$ divides $G_{X^EY^E}^{M^E}-1$.}

\vskip0.2in

Let $\mathscr{L}(\overline{A})=\bigoplus\limits_{{N\in\overline{A}-{\rm
ind}}}\mathbf{Z} u_{[N]}$ be the free Abel group with basis the
set of isomorphism classes determined by indecomposable
$\overline{A}$-modules.

\vskip 0.2in

{\bf Theorem 5.2.}\ {\it $\mathscr{L}(\overline{A})$ is the Lie
subalgebra of ${\mathscr H}(\overline{A})_1$.}

\vskip 0.1in

{\bf Proof.}\  Assume that$X$ and $Y$ are indecomposable
$\overline{A}$-modules such that $X\not\simeq Y$. For any
$\overline{A}$-module $M$, according to Theorem 4.8, we know that
the Hall polynomials $g_{XY}^M$ and $g_{YX}^M\in\mathbf{Z}[x]$
exist, satisfying for any $E\in\Omega_{M\oplus X\oplus Y}$,
$g_{XY}^M(|E|)=G_{X^EY^E}^{M^E}$ and
$g_{YX}^M(|E|)=G_{Y^EX^E}^{M^E}$. By Lemma 5.1, in ${\mathscr
H}(\overline{A})_1$,
$$
u_{[X]}\cdot u_{[Y]}=\sum_{Z\in\overline{A}-{\rm ind}}g_{XY}^Z(1)u_{[Z]}+u_{[X\oplus Y]},
$$
$$
u_{[Y]}\cdot u_{[X]}=\sum_{Z\in\overline{A}-{\rm
ind}}g_{YX}^Z(1)u_{[Z]}+u_{[X\oplus Y]},
$$
therefore
$$
[u_{[X]},u_{[Y]}]=\sum_{Z\in\overline{A}-{\rm
ind}}(g_{XY}^Z(1)-g_{YX}^Z(1))u_{[Z]}\in
\mathscr{L}(\overline{A}),
$$
so we have that
$\mathscr{L}(\overline{A})$ is the Lie subalgebra of
$\mathscr {H}(\overline{A})_1$.     $\hfill\Box$

\vskip 0.2in

Let $\mathscr{L'}(\overline{A})$ be the Lie subalgebra of
$\mathscr{L}(\overline{A})$ generated by the simple
$\overline{A}$-modules. According to [R5] and by using PBW-basis
Theorem, we have the following.

\vskip 0.2in

{ \bf Proposition 5.3. }\ {\it
$\mathscr{C}(\overline{A})_1\otimes_{\mathbf{Z}}\mathbf{Q}$ is the
universal enveloping algebra of
$\mathscr{L'}(\overline{A})\otimes_{\mathbf{Z}}\mathbf{Q}$.}

\vskip 0.2in

Let $\mathscr{L}_0(\overline{A})$ be
the Lie subalgebra of $\mathscr{L'}(\overline{A})$ generated by
$S_1,\cdots,S_n$ and $\mathscr{L}_1(\overline{A})$ the Lie
subalgebra of $\mathscr{L'}(\overline{A})$ generated by
$S_{1'},\cdots,S_{n'}$ respectively. Then by [Rie] we have
$\mathscr{L}_0(\overline{A})\cong\mathscr{L}_1(\overline{A})$ as Lie subalgebras,
which is also isomorphic to the positive part of the corresponding affine Kac-Moody algebra
of type $A$.

\vskip 0.2in

We denote by $\Sigma^{{\rm PI}}_1$ the set of indecomposable projective-injective
$\overline{A}$-modules which are predecessors of $\Sigma_1$,
and by $\Sigma^{{\rm PI}}_2$ the set of indecomposable projective-injective
$\overline{A}$-modules which are successors of $\Sigma_1$.
Note that $\Sigma_1 < \Sigma^{{\rm PI}}_2< \Sigma_2$.

\vskip 0.2in

Let $\Xi_{01} =\{ \ X\in \overline{A}-{\rm ind}\ |\  \Sigma^{{\rm PI}}_1\leq X
\leq \Sigma^{{\rm PI}}_2  \}$. For any $M\in \Xi_{01}$ which is not projective-injective,
reading from the Auslander-Reiten quiver of $\Gamma_{\overline{A}}$, we know that
$\Sigma_1\leq X \leq \tau^{-2}\Sigma_2$. We denote by $\mathscr{L}(\Xi_{01})$ the
free subgroup $\bigoplus\limits_{{N\in\Xi_{01}}}\mathbf{Z} u_{[N]}$ of $\mathscr{L}(\overline{A})$
and by $\mathscr{H}(\Xi_{01})_1$ the subalgebra of $\mathscr{H}(\overline{A})_1$ generated by
indecomposable $\overline{A}$-modules in $\Xi_{01}$.
Let $\mathscr{L}_{01}(\overline{A})= \mathscr{L}(\Xi_{01})\bigcap \mathscr{C}(\overline{A})_1$ and
$\mathscr{C}(\overline{A})_{01}= \mathscr{H}(\Xi_{01})_1\bigcap \mathscr{C}(\overline{A})_1$.

\vskip 0.2in

{ \bf Theorem 5.4. }\ {\it {\rm (1)}\ \
$\mathscr{L}_{01}(\overline{A})$ is a Lie subalgebra of
$\mathscr{L'}(\overline{A})$. In particular,
$\mathscr{L'}(\overline{A})= \mathscr{L}_0(\overline{A})\oplus
\mathscr{L}_{01}(\overline{A})\oplus \mathscr{L}_1(\overline{A})$.}

\vskip 0.1in

{\it {\rm (2)}\ $\mathscr{C}(\overline{A})_{01}\otimes_{\mathbf{Z}}
\mathbf{Q}$ is the universal enveloping algebra of
$\mathscr{L}_{01}(\overline{A})\otimes_{\mathbf{Z}}\mathbf{Q}$.}

\vskip 0.1in

{\bf Proof:}\ (2) is trivial, so we only need to prove (1).
According to the Auslander-Reiten quiver $\Gamma_{\overline{A}}$, we
know that ${\rm add}\ \Xi_{01}$ is closed under extensions, thus
$\mathscr{L}(\Xi_{01})$ is a Lie subalgebra of
$\mathscr{L}(\overline{A})$, hence $\mathscr{L}_{01}(\overline{A})$
is a Lie subalgebra of $\mathscr{L'}(\overline{A})$. The proof is
completed.         $\hfill\Box$

\vskip 0.2in

Let $\delta=(a_1,\cdots,a_n)$ be the minimal positive imaginary root of $A$
and $m=\sum\limits_{i=1}^na_i$.

\vskip 0.2in

{ \bf Theorem 5.5. }\ {\it  Let $M$ be an indecomposable
$\overline{A}$-module and $l(M)$ be the length of $M$. Assume that
$m$ cannot divide $l(M)$, i.e., $m\not | l(M)$, then $u_{[M]}$
belongs to
$\mathscr{C}(\overline{A})_1\otimes_{\mathbf{Z}}\mathbf{Q}$.}

\vskip 0.1in

{\bf Proof:}\  First we assume that $M$ belongs to one of
$\mathscr{P}_0$, $\mathscr{X}_0$, $\mathscr{X}_1$, $\mathscr{I}_1$.
If $M$ is not a projective-injective $\overline{A}$-module, then
$u_{[M]}\in
\mathscr{C}(\overline{A})_1\otimes_{\mathbf{Z}}\mathbf{Q}$ follows
by Theorem 3.2 in [SZ3].

If $M$ is projective-injective, then according to the proof of Theorem 3.3,
we know that $u_{[M]}\in \mathscr{C}(\overline{A})_1\otimes_{\mathbf{Z}}\mathbf{Q}$.

Finally, we assume that $M$  belongs to one of  $\mathscr{R}_0$,  $\mathscr{R}_{01}$,
$\mathscr{R}_1$, that is $M$ belongs to one tube. Note that $m$ cannot divide $l(M)$, it follows that
$M$ must belong to non-homogenous tube, then by Corollary 3.1 in [SZ3], we know that $u_{[M]}$ belongs to
$\mathscr{C}(\overline{A})_1\otimes_{\mathbf{Z}}\mathbf{Q}$. The proof is finished.  $\hfill\Box$

\vskip 0.2in

{\bf Remark.} \  According to Theorem 3.3 and Theorem 5.5, there is
a big difference between the indecomposable $\overline{A}$-modules
which belong to $\mathscr{C}(\overline{A})$ and those which belong
to $\mathscr{C}(\overline{A})_1$, see the following example.

\vskip 0.2in

{\bf Example 5.6.}\  Let $\overline{A}$ be the duplicated tame
quiver algebra of type $\widetilde{D}_4$. That is,
$\overline{A}=\overline{k \widetilde{D_4}}/I$ as in Example 3.4.
According to Theorem 5.5, $u_{[M]}$ belongs to
$\mathscr{C}(\overline{A})_1$ if and only if $M$ belongs to
$\mathscr{P}_0$, $\mathscr{X}_0$, $\mathscr{X}_1$, $\mathscr{I}_1$
or to $\mathscr{R}_0$, $\mathscr{R}_{01}$, $\mathscr{R}_1$ with
$6\not |l(M)$.

\vskip 0.2in

The following example indicates that the converse of Theorem 5.5
does not hold.

\vskip 0.2in

{\bf Example 5.7.} \ Let $K$ be the Kronecker algebra and $\overline{K}$ be the
duplicated algebra of $K$. We may assume that $\overline{K}= kQ_{\overline{K}}/I$ with
$Q_{\overline{K}}: 1\leftleftarrows 2\leftleftarrows
1'\leftleftarrows 2'$. The indecomposable projective-injective
$\overline{A}$-modules are
$P_{1'}= \begin{array}{c}1'\\[-2ex]22\\[-2ex]1\end{array}$ and
$P_{2'}= \begin{array}{c}2'\\[-2ex]1'1'\\[-2ex]2\end{array}$
which are represented by their Loewy series.

\vskip 0.1in

Let $\mathscr{C}(\overline{K})_1$ be the degenerated composition
algebra generated by simple $\overline{K}$-modules
$S_1, S_2, S_{1'}, S_{2'}$. Note that $\delta= (1,1)$ is the minimal positive imaginary root of $K$
and $m=2$ in this case. $l(P_{1'})=4$ and
$u_{[P_{1'}]}\in\mathscr{C}(\overline{A})_1\otimes_{\mathbf{Z}}\mathbf{Q}$ since
$u_{[P_{1'}]}= [u_{[S_{1'}]}, [u_{[S_1]}, u_{[S_2]}u_{[S_2]}]]$.

\vskip 0.1in

Let $M$ be an indecomposable $\overline{K}$-module.
If $l(M)$, the length of $M$, is a positive even number,
then $u_{[M]}\in\mathscr{C}(\overline{A})_1\otimes_{\mathbf{Z}}\mathbf{Q}$ if and only if
$M$ is $P_{1'}$ or $P_{2'}$ since otherwise $M$ belongs to one of homogeneous tubes and in this case
$u_{[M]}$ dose not belong to $\mathscr{C}(\overline{A})_1\otimes_{\mathbf{Z}}\mathbf{Q}$.

\vskip 0.1in

If $l(M)$ is a positive odd number, then from the Auslander-Reiten
quiver of $\Gamma_{\overline{K}}$, we know that $M$ belongs to one
of components $\mathscr{P}_0$, $\mathscr{X}_0$, $\mathscr{X}_1$,
$\mathscr{I}_1$. By using Theorem 3.3, we have that
$u_{[M]}\in\mathscr{C}(\overline{A})_1\otimes_{\mathbf{Z}}\mathbf{Q}$.

\vskip 0.2in

{\bf Acknowledgments.}\ The authors would like to thank Wenxu Ge and Hongbo Lv for many useful discussions.

\begin{description}

\item{[ABST1]}\ I.Assem, T.Br$\ddot{\rm u}$stle, R.Schiffer,
G.Todorov, Cluster categories and duplicated algebras. {\it J.
Algebra}, 305(2006), 548-561.

\item{[ABST2]}\ I.Assem, T.Br$\ddot{\rm u}$stle, R.Schiffer,
G.Todorov,  $m$-cluster categories and $m$-replicated algebras. {\it
Journal of pure and applied Algebra} 212(2008), 884-901.

\item{[ARS]}\ M.Auslander, I.Reiten, S.O.Smal$\phi$, \ Representation
Theory of Artin Algebras. Cambridge Univ. Press, 1995.

\item{[CD]}\ J.Chen, B.Deng, Fundamental relations in Ringel-Hall algebras.
 J.Algebra, 320(2008), 1133 -1149.

\item{[DR]}\  V.Dlab, C.M.Ringel,   Indecomposable representations of graphs and
algebras. Mem.Amer.Math.Soc., 173 (1976)

\item{[G]}\  J.A.Green, Hall algebras,hereditary algebras and quantum groups.
Inventiones Math., 120(1995). 361-377.

\item{[HRS]}\ D.Happel, I.Reiten, S.O.Smal$\phi$,  Tilting in abelian
   categories and quasitilted algebras, Mem.Amer.Math.Soc., 575, 1996.

\item{[La]}\  T.Y.Lam, First Course in Noncommutative
Rings.  Graduate Texts in Mathematics, Vol. 131, Springer-Verlag,
New York, 1991.

\item{[Lu]}\  G.Lusztig, Intoduction to quantum groups. Birkh\H{a}user.
   Boston.1993.

     \item{[P]}\  L.Peng, Some Hall polynomials for
    representation-finite trivial extension algebras.  J.Algebra,
    197(1997),1-13.

\item{[PX]}\  L.Peng, J.Xiao, Triangulated categories and Kac-Moody algebras.
   Invent. math., 140(2000), 563-603.

   \item{[R1]}\  C.M.Ringel,  Tame algebras and integral quadratic forms. Springer
   Lecture Notes in Math.,  (1099), 1984

    \item{[R2]}\ C.M.Ringel, Hall algebra.  In: Topics in Algebras.
   Banach Centre Publ. Warszawa. 1990, \ 26, \ 433-447.

    \item{[R3]}\ C.M.Ringel, Hall algebras and quantum groups.
   Inventiones Math.,  101 (1990), 583-592

    \item{[R4]}\ C.M.Ringel, Green's theorem on Hall algebras.  Canad.Math.Soc.
   Conf.Proc. Vol 19, pp:185-245.

   \item{[R5]}\ C.M.Ringel, Lie algebra arising in representation
theory, London Math Soc., 168, Cambridge University Press, 1992,
284-291.

\item{[Rie]}\ C.H.Riedtmann, Lie algebras generated by indecomposables.
   J.Algebra, 170(1994), 526-546.

   \item{[PZ1]}\ P.Zhang, Composition algebras of affine type.
   J.Algebra.206(1998),505-540.

   \item{[PZ2]}\ P.Zhang, Rigids as iterated skew commutators of simples.
   Algebr Represent Theory, 9(2006), 539-555.

    \item{[ZZ]}\ P.Zhang, S.Zhang, Indecomposables as elements in affine
   composition algebras.  J.Algebra, 210(1998), 614-629.

    \item{[SZ1]}\ S.Zhang, The Hall polynomials for tame quiver algebras. J.Algebra, 239(2001), 606-614.

    \item{[SZ2]}\ S.Zhang, Triangular decomposition of  tame non-simple-laced
   composition algebras. J.Algebra, 241(2001),548-577.

   \item{[SZ3]}\ S.Zhang, Lie algebra determined by tame hereditary algebras (In chinese).
   Chinese Ann. Math. 21(5)(A)(2000), 619-621.

\end{description}

\end{document}